\documentclass{article}
\usepackage{amsthm}
\usepackage{epsfig,amsmath,amssymb}

\newtheorem{teo}{Theorem}[section]
\numberwithin{equation}{section}

\newtheorem{cor}[teo]{Corollary}
\newtheorem{pro}[teo]{Proposition}
\newtheorem{lem}[teo]{Lemma}

\theoremstyle{definition}
\newtheorem{defi}[teo]{Definition}
\newtheorem{obs}[teo]{Remark}

\title{Right Coideal Subalgebras of the Quantum Borel Algebra of type $G_2$\thanks{This paper was written during a research period at the UNAM FES-C, Mexico, with the support of CNPq-Brazil. It will be part of the author's PhD thesis.}}

\author{B\'arbara Pogorelsky\\\footnotesize{barbara.pogorelsky@ufrgs.br}\\\footnotesize{Instituto de Matem\'atica, Universidade Federal do Rio Grande do Sul}\\\footnotesize{Av. Bento Gon\c calves 9500, Porto Alegre, RS, 91509-900, Brazil}}

\begin{document}

\maketitle

\begin{abstract}
\noindent In this paper we describe the right coideal subalgebras
containing all group-like elements of the multiparameter quantum
group $U_q^+(\mathfrak{g})$, where $\mathfrak{g}$ is a simple Lie
algebra of type $G_2$, while the main parameter of quantization
$q$ is not a root of 1. If the multiplicative order $t$ of $q$ is
finite, $t>4$, $t\neq 6$, then the same classification remains
valid for homogeneous right coideal subalgebras of the positive
part $u_q^+(\mathfrak{g})$ of the multiparameter version of the
small Lusztig quantum group.
\end{abstract}

{\scriptsize{\it Keywords:} Quantum groups, Hopf algebras, Coideal
subalgebras, PBW-generators}

{\scriptsize{\it Mathematics Subject Classification:} 16W30, 17B37}

\section{Introduction}

\quad Comodule algebras over a Hopf algebra naturally arise in the
Galois theory of Hopf algebra actions as the Galois objects; see
A. Masuoka and T. Yanai \cite{MY}, A. Milinski \cite{Mil}, S.
Westreich and T. Yanai \cite{WY} and T. Yanai \cite{Y, Y2}. In
particular the Galois correspondence theorem for the actions on a
free algebra sets up a one to one correspondence between all right
coideal subalgebras and all intermediate free subalgebras; see
V.O. Ferreira, L.S.I. Murakami, and A. Paques \cite{FMP}. At the
same time, the notion of one-sided coideal subalgebras appears to
be of fundamental importance in the theory of quantum groups: a
survey by G. Letzter \cite{Let} provides an overview of the use of
one-sided coideal subalgebras in constructing quantum symmetric
pairs, in forming quantum Harish-Chandra modules and in producing
quantum symmetric spaces.

Recently V. K. Kharchenko and A. V. Lara Sagah\'on \cite{Kh6},
using a PBW-basis construction method \cite{Kh5}, offered a
complete classification of right coideal subalgebras that contain
the coradical $\textbf{k}[G]$ for the quantum group
$U_q(\mathfrak{sl}_{n+1})$. As a consequence they determined that
the quantum Borel algebra $U_q^+(\mathfrak{g})$,
$\mathfrak{g}=\mathfrak{sl}_{n+1}$, contains $(n+1)!$ different
right coideal subalgebras that include the coradical. If
$\mathfrak{g}=\mathfrak{so}_{2n+1}$ is a simple Lie algebra of
type $B_n$ then $U_q^+(\mathfrak{g})$ has $(2n)!!$ right coideal
subalgebras that include the coradical \cite{Kh8}. In both cases
the number coincides with the order of the Weyl group defined by
the Lie algebra $\mathfrak{g}$. This provides enough reason to
conjecture that for arbitrary simple finite dimensional Lie
algebras $\mathfrak{g}$ the number of right coideal subalgebras in
$U_q^+(\mathfrak{g})$ that include the coradical coincides with
the order of the Weyl group related to $\mathfrak{g}$, see
\cite{Kh8}.

In this paper by means of the same PBW-basis construction method
we prove this conjecture for the Lie algebra $\mathfrak{g}$ of
type $G_2$. More precisely, we prove the following theorem.

\begin{teo}\label{lattice}
If $q$ is not a root of 1, the lattice of right coideal
subalgebras containing \emph{$\textbf{k}[G]$} of
$U_q^+(\mathfrak{g})$ is given in the Figure 1. If $q$ has
multiplicative order $t>4$, $t \neq 6$, the same figure is the
lattice of homogeneous right coideal subalgebras containing
\emph{$\textbf{k}[G]$} of $u_q^+(\mathfrak{g})$.

\begin{figure}[h]
\begin{center}
\unitlength 1mm
\begin{picture}(60.00,81.25)(0,0)

\put(35.00,81.25){\makebox(0,0)[cc]{$U_q^+(\mathfrak{g})$}}
\put(20.00,71.25){\makebox(0,0)[cc]{}}
\put(60.00,61.25){\makebox(0,0)[bc]{$\langle[x_1,x_2]\rangle$}}
\put(0.00,5.00){\makebox(0,0)[cc]{}}
\put(60.00,51.25){\makebox(0,0)[bc]{$\langle[[x_1,x_2],[[x_1,x_2],x_2]]\rangle$}}
\put(60.00,41.25){\makebox(0,0)[bc]{$\langle[[x_1,x_2],x_2]\rangle$}}
\put(60.00,31.25){\makebox(0,0)[bc]{$\langle[[[x_1,x_2],x_2],x_2]\rangle$}}
\put(60.00,21.25){\makebox(0,0)[bc]{$\langle x_{2}\rangle$}}
\put(35.00,5.63){\makebox(0,0)[cc]{$\textbf{k}[G]$}}
\put(55.00,0.00){\makebox(0,0)[cc]{}}
\linethickness{0.15mm}
\put(60.00,55.63){\line(0,1){4.38}}
\put(10.00,51.25){\makebox(0,0)[bc]{$\langle[x_{2},[x_{2},[x_{2},x_{1}]]]\rangle$}}
\put(10.00,41.25){\makebox(0,0)[bc]{$\langle[x_{2},[x_{2},x_{1}]]\rangle$}}
\put(10.00,31.25){\makebox(0,0)[bc]{$\langle[x_{2},x_{1}]\rangle$}}
\put(10.00,21.25){\makebox(0,0)[bc]{$\langle x_{1}\rangle$}}
\linethickness{0.15mm}
\put(10.00,61.25){\makebox(0,0)[bc]{$\langle[[x_{1},x_{2}],[x_{2},[x_{2},x_{1}]]]\rangle$}}
\linethickness{0.15mm}
\put(60.00,45.63){\line(0,1){4.38}}
\linethickness{0.15mm}
\put(60.00,35.63){\line(0,1){4.38}}
\linethickness{0.15mm}
\put(60.00,25.63){\line(0,1){4.38}}
\linethickness{0.15mm}
\put(10.00,55.63){\line(0,1){4.38}}
\linethickness{0.15mm}
\put(10.00,45.63){\line(0,1){4.38}}
\linethickness{0.15mm}
\put(10.00,35.63){\line(0,1){4.38}}
\linethickness{0.15mm}
\put(10.00,25.63){\line(0,1){4.38}}
\linethickness{0.15mm}
\multiput(40.00,80.00)(0.18,-0.12){109}{\line(1,0){0.18}}
\linethickness{0.15mm}
\multiput(10.00,20.00)(0.18,-0.12){109}{\line(1,0){0.18}}
\linethickness{0.15mm}
\multiput(10.00,66.25)(0.17,0.12){115}{\line(1,0){0.17}}
\linethickness{0.15mm}
\multiput(40.00,6.25)(0.17,0.12){115}{\line(1,0){0.17}}
\end{picture}
\caption{Lattice of Right Coideal Subalgebras}
\end{center}
\end{figure}
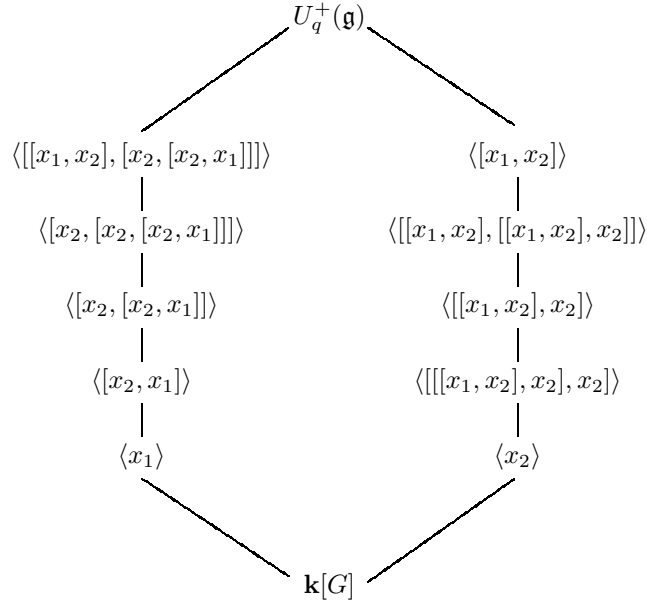
\end{teo}

Here the distinguished element is a generator of the right coideal
subalgebra (in particular each right coideal subalgebra is
generated over the coradical by a single element).

In the second section, following to \cite{Kh6}, we introduce main
concepts and general results that are of use for further
considerations. In the third section we consider the algebra
$U_q^+(\mathfrak{g})$ (respectively, $u_q^+(\mathfrak{g})$) as a
character Hopf algebra \cite{Kh4} in order to find a PBW-basis in
the explicit form (Theorem \ref{base}). The results of this section are very similar to the results of
the section 4.2 of the paper \cite{A} by I. Angiono dedicated to classification of
finite dimensional Nichols algebras (quantum symmetric algebras) over
algebraically closed fields of characteristic zero. In the fourth section,
following \cite{Kh5}, we transform the found PBW-basis up to a
PBW-basis of a given right coideal subalgebra. In this way we may
find all possible PBW-bases for right coideal subalgebras
containing $\textbf{k}[G]$. This provides the required
classification (Theorem \ref{lattice}).

\section{Preliminaries}

\begin{defi}
Let $S$ be an algebra over a field $\textbf{k}$ and $A$ its
subalgebra with a fixed basis $\{a_j|j\in J\}$. A linearly ordered
subset $W \subseteq S$ is said to be a set of
\textit{PBW-generators of $S$ over $A$} if there exists a function
$h:W\rightarrow\mathbb{Z}^+\cup\infty$, called the height
function, such that the set of all products
\begin{equation}\label{pbw}
a_jw_1^{n_1}w_2^{n_2}\ldots w_k^{n_k},
\end{equation}
where $j \in J$, $w_1<w_2<\ldots <w_k \in W$, $n_i<h(w_i)$, $1\leq
i\leq k$ is a basis of $S$. The value $h(w)$ is referred to as the
\textit{height} of $w$ in $W$. If $A=\textbf{k}$ is the ground
field, then we shall call $W$ simply as a set of PBW-generators of
$S$.
\end{defi}

\begin{defi}
Let $W$ be a set of PBW-generators of $S$ over a subalgebra $A$.
Suppose that the set of all words in $W$ as a free monoid has its
own order $\prec$ (that is, $a \prec b$ implies $cad \prec cbd$
for all words $a,b,c,d \in W$). A \textit{leading word} of $s \in
S$ is the maximal word $m=w_1^{n_1}w_2^{n_2}\ldots w_k^{n_k}$ that
appears in the decomposition of $s$ in the basis (\ref{pbw}). A
\textit{leading term} of $s$ is the sum $am$ of all terms
$\alpha_ia_im$, $\alpha_i \in \textbf{k}$, that appear in the decomposition of $s$ in the
basis (\ref{pbw}), where $m$ is the leading word of $s$.
\end{defi}

\begin{defi}
A Hopf algebra $H$ is said a \textit{character Hopf algebra} if
the group $G$ of all group-like elements is commutative and $H$ is
generated over $\textbf{k}[G]$ by skew primitive semi-invariants
$a_i, i \in I$: $$\Delta(a_i)=a_i\otimes 1+g_i\otimes a_i,\quad
g^{-1}a_ig=\chi^i(g)a_i,\quad g, g_i \in G,$$ where $\chi^i$, $i
\in I$, are characters of the group $G$.
\end{defi}

Let us associate a quantum variable $x_i$ to $a_i$. For each word
$u$ in $X=\{x_i|i\in I\}$ we denote by $g_u$ an element of $G$
that appears from $u$ by replacing each $x_i$ with $g_i$. In the
same way we denote by $\chi^u$ a character that appears from $u$
by replacing each $x_i$ with $\chi^i$. We define a bilinear skew
commutator on homogeneous linear combinations of words by the
formula
\begin{equation}\label{brac}
[u,v]=uv-\chi^u(g_v)vu,
\end{equation}
where we use the notation $\chi^u(g_v)=p_{uv}=p(u,v)$. These
brackets are related to the product by the following identities
\begin{equation}\label{bracprod1}
[u\cdot v,w]=p_{vw}[u,w]\cdot v+u\cdot [v,w],
\end{equation}
\begin{equation}\label{bracprod2}
[u,v\cdot w]=[u,v]\cdot w+p_{uv}v\cdot [u,w].
\end{equation}

The group $G$ acts on the free algebra $\textbf{k}\langle
X\rangle$ by $g^{-1}ug=\chi(g)u$, where $u$ is an arbitrary
monomial in $X$. The skew group algebra $G\langle X\rangle$ has
the natural Hopf algebra structure
$$\Delta(x_i)=x_i\otimes 1+g_i\otimes x_i,\quad i \in I,\quad
\Delta(g)=g\otimes g.$$ We fix a Hopf algebra homomorphism
$$\xi :G\langle X\rangle\rightarrow H,\quad \xi(x_i)=a_i,\quad
\xi(g)=g,\quad i\in I,\quad g\in G.$$

\begin{defi}
A \textit{constitution} of a word $u$ in $G\cup X$ is a family of
non-negative integers $\{m_x,x\in X\}$ such that $u$ has $m_x$
occurrences of $x$. Certainly almost all $m_x$ in the constitution
are zero.
\end{defi}

Let us fix an arbitrary complete order $<$ on the set $X$, and let
$\Gamma^+$ be the free additive (commutative) monoid generated by
$X$. The monoid $\Gamma^+$ is a completely ordered monoid with
respect to the following order:
\begin{equation}\label{ordemgrau}
m_1x_{i_1}+m_2x_{i_2}+\ldots
+m_kx_{i_k}>m_1'x_{i_1}+m_2'x_{i_2}+\ldots +m_k'x_{i_k}
\end{equation}
if the first from the left nonzero number in
$(m_1-m_1',m_2-m_2',\ldots,m_k-m_k')$ is positive, where
$x_{i_1}>x_{i_2}>\ldots>x_{i_k}$ in $X$. We associate a formal
degree $D(u)=\sum_{x \in X}m_xx \in \Gamma^+$ to a word $u$ in
$G\cup X$, where $\{m_x|x\in X\}$ is the constitution of $u$.
Respectively, if $f=\sum \alpha_iu_i \in G\langle X\rangle$,
$0\neq \alpha_i \in \textbf{k}$ then
\begin{equation}\label{grau}
D(f)=max_i\{D(u_i)\}. \end{equation}
On the set of all words in $X$
we fix the lexicographical order with the priority from the left
to the right, where a proper beginning of a word is considered to
be greater than the word itself.

\begin{defi}
A non-empty word $u$ is called a \textit{standard word} (or
\textit{Lyndon word}, or \textit{Lyndon-Shirshov word}) if $vw>wv$
for each decomposition $u=vw$ with non-empty $v,w$.
\end{defi}

\begin{defi}
A \textit{non-associative word} is a word where brackets $[,]$ are
somehow arranged to show how multiplication applies.
\end{defi}

If $[u]$ denotes a non-associative word, then by $u$ we denote an
associative word obtained from $[u]$ by removing the brackets. Of
course, $[u]$ is not uniquely defined by $u$ in general.

\begin{defi}
The set of \textit{standard non-associative words} is the biggest
set $SL$ that contains all variables $x_i$ and satisfies the
following properties:
\begin{enumerate}
    \item If $[u]=[[v],[w]] \in SL$ then $[v],[w] \in SL$, and
    $v>w$ are standard.
    \item If $[u]=[[[v_1],[v_2]],[w]] \in SL$ then $v_2\leq w$.
\end{enumerate}
\end{defi}

By the Shirshov's Theorem, every standard word has only one
alignment of brackets such that the defined non-associative word
is standard. In order to find this alignment we use the following
procedure: the factors $v,w$ of the non-associative decomposition
$[u]=[[v],[w]]$ are standard words such that $u=vw$ and $v$ has
the minimal length (see \cite{sh}).

\begin{defi}
A \textit{super-letter} is a polynomial that equals a
non-associative standard word where the brackets mean
(\ref{brac}). A \textit{super-word} is a word in super-letters.
\end{defi}

Using Shirshov's Theorem, every standard word $u$ defines only one
super-letter that will be denoted by $[u]$. The order on the
super-letters is defined in the natural way:
$[u]>[v]\Leftrightarrow u>v$.

Since quantum Borel algebras $U_q^+(\mathfrak{g})$ and
$u_q^+(\mathfrak{g})$ are homogeneous in each variable, in what
follows we suppose that $H$ is a $\Gamma^+$-graded character Hopf
algebra, that is, $H$ is homogeneous in each of the generators
$a_i$.

\begin{defi}\label{harddef}
A super-letter $[u]$ is called \textit{hard in $H$} if its value
in $H$ is not a linear combination of super-words of the same
degree \eqref{grau} in super-letters smaller than $[u]$.
\end{defi}

\begin{pro}\label{hardeq}
\emph{(\cite[Corollary 2]{Kh2})} A super-letter $[u]$ is hard in
$H$ if and only if the value in $H$ of the standard word $u$ is
not a linear combination of values of smaller words of the same
degree (\ref{grau}).
\end{pro}

\begin{pro}\label{4.8}
\emph{(\cite[Lemma 4.8]{Kh4})} Let $B$ be a set of super-letters
containing $x_1,\ldots,x_n$. If each pair $[u],[v] \in B$, $u>v$
satisfies one of the following conditions
\begin{description}
    \item[1)] $[[u],[v]]$ is not a standard non-associative word;
    \item[2)] the super letter $[[u],[v]]$ is not hard in $H$;
    \item[3)] $[[u],[v]] \in B$;
\end{description}
then the set $B$ includes all hard in $H$ super-letters.
\end{pro}

\begin{defi}\label{altura}
We say that the \textit{height} of a hard in $H$ super-letter
$[u]$ equals $h=h([u])$ if $h$ is the smallest number such that
\begin{enumerate}
    \item $p_{uu}$ is a primitive $t$-th root of $1$ and
    either $h=t$ or $h=tl^r$, where $l=char(\textbf{k})$,
    \item the value of $[u]^h$ in $H$ is a linear combination of
    super-words of the same degree (\ref{grau}) in super-letters smaller than
    $[u]$.
\end{enumerate}
If there exists no such number then the height equals infinity.
\end{defi}

\begin{teo}\label{hard}
\emph{(\cite[Theorem 2]{Kh2})} The values of all hard in $H$
super-letters with the above defined height function form a set of
PBW-generators for $H$ over $\emph{\textbf{k}[G]}$.
\end{teo}

According to \cite[Theorem 1.1]{Kh5}, every right coideal
subalgebra $\textbf{U}$ that contains all group-like elements has
a PBW-basis over $\textbf{k}[G]$ which can be extended up to a
PBW-basis of $H$. The PBW-generators $T$ for $\textbf{U}$ can be
obtained from the PBW-basis of $H$ given in Theorem \ref{hard} in
the following way.

Suppose that for a given hard super-letter $[u]$ there exists an
element $c \in \textbf{U}$
\begin{equation}\label{geradores}
c=[u]^s+\sum \alpha_iW_i
\end{equation}
where $W_i$ are the basis super-words starting with super-letters
smaller than $[u]$, $D(W_i)=sD(u)$. We fix one of the elements
with minimal $s$ and denote it by $c_u$. Thus, for every
super-letter $[u]$ hard in $H$ we have at most one element $c_u$.
We define the height function by the following lemma.

\begin{lem}\label{s1t}
\emph{(\cite[Lemma 4.3]{Kh5})} In the representation
\eqref{geradores} of the chosen element $c_u$, either $s=1$ or
$p(u,u)$ is a primitive $t$-th root of $1$ and $s=t$, or (in the
case of positive characteristic) $s=t(char \emph{\textbf{k}})^r$.
\end{lem}

If the height of $[u]$ in $H$ is infinite, then the height of
$c_u$ in $\textbf{U}$ is defined to be infinite as well. If the
height of $[u]$ in $H$ equals $t$ and $p(u,u)$ is a primitive
$t$-th root of $1$, then, due to the above lemma, $s=1$ (note that
in the representation \eqref{geradores} the number $s$ is less
than the height of $[u]$). In this case, the height of $c_u$ in
$\textbf{U}$ is supposed to be $t$ as well. If the characteristic
$l$ is positive and the height of $[u]$ in $H$ equals $tl^r$, then
we define the height of $c_u$ in $\textbf{U}$ to be equal to
$tl^r/s$ (thus, in characteristic zero the height of $c_u$ in
$\textbf{U}$ always equals the height of $[u]$ in $H$).

\begin{pro}\label{pbwU}
\emph{(\cite[Proposition 4.4]{Kh5})} The set of all chosen $c_u$
with the above defined height function forms a set of
PBW-generators for \emph{$\textbf{U}$} over
\emph{$\textbf{k}[G]$}.
\end{pro}

\begin{defi}\label{U+} (see, for example, \cite[section 3]{Kh6})
Let $C=\parallel a_{ij} \parallel$ be a generalized Cartan matrix
symmetrizable by $D=diag(d_1,\ldots,d_n)$, $d_ia_{ij}=d_ja_{ji}$.
Denote by $\mathfrak{g}$ a Kac-Moody algebra defined by $C$ (see
\cite{kac}). Suppose that the quantification parameters
$p_{ij}=p(x_i,x_j)=\chi^i(g_j)$ are related by
\begin{equation}\label{pis}
p_{ii}=q^{d_i},\quad p_{ij}p_{ji}=q^{d_ia_{ij}},\quad 1\leq i,j
\leq n.
\end{equation}
The \textit{multiparameter quantization $U^+_q(\mathfrak{g})$} of
the Borel subalgebra $\mathfrak{g^+}$ is a character Hopf algebra generated by $x_1,\ldots,x_n,g_1,\ldots,g_n$ and
defined by Serre relations with the skew brackets (\ref{brac}) in
place of the Lie operation:
\begin{equation}\label{serre}
[[\ldots[[x_i,x_j],x_j],\ldots],x_j]=0,\quad 1\leq i\neq j\leq n,
\end{equation}
where $x_j$ appears $1-a_{ji}$ times.
\end{defi}

By \cite[Theorem 6.1]{Kh1} the left hand sides of these relations
are skew-primitive elements in $G\langle X \rangle$. Therefore,
the ideal generated by these elements is a Hopf ideal, while
$U^+_q(\mathfrak{g})$ indeed has a natural character Hopf algebra
structure.

\begin{defi}(see, for example, \cite[section 3]{Kh6})
If the multiplicative order $t$ of $q$ is finite, then we define
$u^+_q(\mathfrak{g})$ as $G\langle X\rangle/\Lambda$, where
$\Lambda$ is the biggest Hopf ideal in $G\langle X\rangle^{(2)}$,
which is the set (an ideal) of noncommutative polynomials without
free and linear terms. From \cite[Lemma 2.2]{Kh7} this is a
$\Gamma^+$-homogeneous ideal. Certainly $\Lambda$ contains all
skew-primitive elements of $G\langle X\rangle^{(2)}$ (each one of
them generates a Hopf ideal). Hence, by \cite[Theorem 6.1]{Kh1},
relations (\ref{serre}) are still valid in $u^+_q(\mathfrak{g})$.
\end{defi}

We notice that the subalgebra $A$ generated by $x_1,\ldots,x_n$ over $\textbf{k}$ in
$U_q^+(\mathfrak{g})$ is a Nichols algebra of Cartan type if $q$ is not a root of $1$, as described in \cite{AS}. Analogously, if $q^t=1$ for an integer $t$, the same thing is valid for $A \subseteq u_q^+(\mathfrak{g})$. This is particularly useful since in \cite{A} there are many results for the Nichols algebra $A$, although in this paper {\bf  k} is an  algebraically
closed field of characteristic zero. However, if $q$ is a root of 1, then
the subalgebra generated by $x_1,\ldots,x_n$ in $U_q({\mathfrak g})$ is not a
Nichols algebra.

\begin{defi}\label{defder}
The subalgebra $A$ generated by $x_1,\ldots,x_n$ over $\textbf{k}$ in
$U_q^+(\mathfrak{g})$ (respectively, $u_q^+(\mathfrak{g})$) has a
\textit{differential calculus} defined by
\begin{equation*}
\partial_i(x_j)=\delta_i^j,\quad
\partial_i(uv)=\partial_i(u)v+p(u,x_i)u\partial_i(v),
\end{equation*}
for $x_i \in X$.
\end{defi}

\begin{lem}\label{derivadasq1}
\emph{(\cite[Lemma 2.10]{Kh8})} Let $u \in$
\emph{$\textbf{k}\langle X\rangle$} be a homogeneous in each $x_i$
element. If $p_{uu}$ is a $t$-th primitive root of $1$, then
\begin{equation*}
\partial_i(u^t)=p(u,x_i)^{t-1}\underbrace{[u,[u,
\ldots[u}_{t-1},\partial_i(u)]\ldots]].
\end{equation*}
\end{lem}

\begin{lem}\label{MS}
\emph{(Milinski-Schneider criterion, see \cite{MS})} If a
polynomial $f \in$ \emph{$\textbf{k}\langle X\rangle$} with no
free term is such that $\partial_i(f)=0$ in $u_q^+(\mathfrak{g})$
for every $x_i \in X$, then $f=0$ in $u_q^+(\mathfrak{g})$.
\end{lem}

\section{Explicit PBW-Generators for Quantizations}

\quad In this section we are going to explicit a set of
PBW-generators for $U_q^+(\mathfrak{g})$ (respectively,
$u_q^+(\mathfrak{g})$, if $q^t=1$ for $t>4$, $t\neq 6$), where $\mathfrak{g}$
is the simple Lie algebra of type $G_2$.

Let us first remember that the algebra $U_q^+(\mathfrak{g})$ is
defined by two generators $x_1,x_2$ and two relations
\begin{equation}\label{rel}
[[[[x_1,x_2],x_2],x_2],x_2]=0,\quad [x_1,[x_1,x_2]]=0,
\end{equation}
where the brackets mean the skew commutator (\ref{brac}).
Relations \eqref{pis} take up the form $p_{11}=p_{22}^3$,
$p_{12}p_{21}=p_{11}^{-1}$, and $p_{22}=q$. In what follows we
shall suppose that $q^6\neq 1$ and $q^4\neq 1$.

Using (\ref{brac}), the relations in (\ref{rel}) can also be
written as
\begin{equation}\label{rel2}
x_1x_2^4+a_1x_2x_1x_2^3+a_2x_2^2x_1x_2^2+a_3x_2^3x_1x_2+a_4x_2^4x_1=0,
\end{equation}
\begin{equation}\label{rel3}
x_1^2x_2+b_1x_1x_2x_1+b_2x_2x_1^2=0,
\end{equation}
where
\begin{equation*}
a_1=-p_{12}p_{22}^{[4]},\quad
a_2=p_{12}^2p_{22}p_{22}^{[3]}(p_{22}^2+1),\quad
a_3=-p_{12}^3p_{22}^3p_{22}^{[4]},\quad a_4=p_{12}^4p_{22}^6,
\end{equation*}
\begin{equation*}
b_1=-p_{12}(1+p_{11}),\quad b_2=p_{12}^2p_{11},
\end{equation*}
and as usual we denote $p^{[n]}=1+p+...+p^{n-1}$ (another usual notation is $[n]_p$).

Let us multiply (\ref{rel2}) from the left by $x_1$, while
(\ref{rel3}) from the right by $x_2^3$. The difference of the
obtained relations provides a new one
\begin{equation}\label{dif1}
(a_1-b_1)x_1x_2x_1x_2^3=-a_2x_1x_2^2x_1x_2^2-a_3x_1x_2^3x_1x_2-a_4x_1x_2^4x_1+b_2x_2x_1^2x_2^3.
\end{equation}

We may multiply (\ref{rel2}) from the left by $(a_1-b_1)x_1x_2$,
while (\ref{dif1}) from the right by $x_2$. Again the difference
gives a new relation
$$\{a_1(a_1-b_1)-a_2\}x_1x_2^2x_1x_2^3=\{a_3-a_2(a_1-b_1)\}x_1x_2^3x_1x_2^2$$
\begin{equation}\label{dif2}
+\{a_4-a_3(a_1-b_1)\}x_1x_2^4x_1x_2-a_4(a_1-b_1)x_1x_2^5x_1-b_2x_2x_1^2x_2^4.
\end{equation}

In the same way, if we multiply (\ref{rel3}) from the right by
$(a_1-b_1)x_1x_2^3$, while (\ref{dif1}) from the left by $x_1$,
the difference of the obtained relations after replacement of all
subwords $x_1^2x_2$ with $-b_1x_1x_2x_1-b_2x_2x_1^2$ defines a new
relation
$$(a_2b_1-\{b_1(a_1-b_1)+b_2\}b_1)x_1x_2x_1x_2x_1x_2^2=\{b_1(a_1-b_1)+b_2\}b_2x_1x_2^2x_1^2x_2^2$$
\begin{equation}\label{dif4}
-a_3b_1x_1x_2x_1x_2^2x_1x_2-a_4b_1x_1x_2x_1x_2^3x_1+W.
\end{equation}
where $W$ is a linear combination of words with the first letter
$x_2$.

Now we are ready to prove the following statement.

\begin{pro}\label{base}
All hard in $U_q^+(\mathfrak{g})$ $($respectively, in
$u_q^+(\mathfrak{g}))$ super-letters are contained in the
following list:
\begin{align}\label{base2}
[A]&=x_1,\notag\\ [B]&=[x_1,x_2],\notag\\
[C]&=[[x_1,x_2],[[x_1,x_2],x_2]],\\ [D]&=[[x_1,x_2],x_2],\notag\\
[E]&=[[[x_1,x_2],x_2],x_2],\notag\\ [F]&=x_2.\notag
\end{align}
\end{pro}

\emph{Proof.} To prove the proposition we'll show that this set
satisfies the conditions of Proposition \ref{4.8}, that is, for
every pair $[X],[Y]$ in this set such that $X
> Y$, the non-associative word $[[X],[Y]]$ either belongs to this
set, or it is not a standard non-associative word, or it defines a
non-hard in $U_q^+(\mathfrak{g})$ super-letter. If one of these 3
possibilities is satisfied for every possible pair, then
Proposition \ref{4.8} proves that all the hard in
$U_q^+(\mathfrak{g})$ super-letters are contained in the list
\eqref{base2}.

We have 15 possibilities. In four cases $[[A],[F]]=[B]$,
$[[B],[D]]=[C]$, $[[B],[F]]=[D]$ and $[[D],[F]]=[E]$. In five more
cases ($[[A],[B]]$, $[[A],[C]]$, $[[A],[D]]$, $[[A],[E]]$ and
$[[E],[F]]$) the non-associative word is not hard by Proposition
\ref{hardeq}, since the associative word obtained by omitting the
brackets contains either subword $x_1x_2^4$ or subword $x_1^2x_2$,
and therefore by the relations (\ref{rel2}) and (\ref{rel3}) it is
a linear combination of smaller words in $U_q^+(\mathfrak{g})$. In
two cases, $[[C],[E]]$ and $[[C],[F]]$, the word is not standard
as a non-associative word. It remains to consider the following
four cases:
\begin{itemize}
    \item $[[D],[E]]$: the relation (\ref{dif2}) shows that the
    word $DE$ is a linear combination of smaller words in $U_q^+(\mathfrak{g})$. Thus, by Proposition \ref{hardeq}, the super-letter is not
    hard.
    \item $[[B],[C]]$: since $BC=(x_1x_2)^3x_2$, the relation
    (\ref{dif4}) shows that the word $BC$ is a linear combination of smaller
    words. Again, by Proposition \ref{hardeq}, the super-letter is not
    hard.
    \item $[[B],[E]]$: in this case $BE=x_1x_2x_1x_2^3$, and we
    may use the relation (\ref{dif1}).
    \item $[[C],[D]]$: let us multiply the relation (\ref{dif4})
    from the right by $(a_1-b_1)x_2$, while the relation
    (\ref{dif1}) from the left by
    $(a_2b_1-\{b_1(a_1-b_1)+b_2\}b_1)x_1x_2$. The leading term of
    the difference equals
    $$\{a_3b_1(a_1-b_1)-a_2(a_2b_1-\{b_1(a_1-b_1)+b_2\}b_1)\}x_1x_2x_1x_2^2x_1x_2^2$$
    $$=-p_{12}^5p_{22}^5(1+p_{22}^3)(1+p_{22}^2)CD.$$
    Therefore $CD$ is also a linear combination of smaller words.
\end{itemize}

Thus, by Proposition \ref{4.8}, the set $\{[A], [B], [C], [D],
[E], [F]\}$ contains all hard in $U_q^+(\mathfrak{g})$
super-letters.

Since $u_q^+(\mathfrak{g})$ is a homomorphic image of
$U_q^+(\mathfrak{g})$, all non-hard in $U_q^+(\mathfrak{g})$
super-letters are non-hard in $u_q^+(\mathfrak{g})$. Hence the
list \eqref{base2} contains all hard in $u_q^+(\mathfrak{g})$
super-letters as well. \hfill $\square$ \vspace{0.5cm}

Here we have to make an important remark. A PBW-basis over $\textbf{k}$ for the Nichols algebra of Cartan type $G_2$ is computed in \cite{A}, and as a consequence we have a PBW-basis over $\textbf{k}[G]$ for $U_q^+(\mathfrak{g})$ provided that $q$ is not a root of 1, and for
$u_q({\mathfrak g})$ otherwise. However, the previous proposition is necessary since to find the PBW-generators for the possible right coideal subalgebras $\textbf{U}$ we need not just a PBW-basis, which is not unique, but the exact PBW-basis constituted by the hard super-letters, as described in Theorem \ref{hard} and Proposition \ref{pbwU}.

Note that, according to Definition \ref{defder}, the subalgebra
$A=\textbf{k}\langle x_1,x_2\rangle$ of $U_q^+(\mathfrak{g})$
(respectively, $u_q^+(\mathfrak{g})$) has a differential calculus
\begin{equation}\label{der}
\partial_i(x_j)=\delta_i^j,\quad \partial_i(uv)=\partial_i(u)\cdot v+p(u,x_i)u\cdot\partial_i(v)
\end{equation}
for $i=1,2$.

Now we notice that, if $\partial_i(u)=0$, then:
\begin{equation}\label{formula}
\partial_i([u,x_i])=\partial_i(u)\cdot x_2+p(u,x_i)u\cdot\partial_i(x_i)-p(u,x_i)(\partial_i(x_i)\cdot u+p(x_i,x_i)x_i\cdot\partial_i(u))=0.
\end{equation}

From \cite[Lemma 2.10]{Kh6} we know that the homogeneous right
coideal subalgebras containing $\textbf{k}[G]$ have the form
$\textbf{U}=U_A\#\textbf{k}[G]$, where $U_A=A\cap\textbf{U}$ is a
differential subalgebra of $A=\textbf{k}\langle x_1,x_2\rangle$
and $G$ is the set of all grouplike elements. So, from knowing the
differential subalgebras, we can describe the homogeneous right
coideal subalgebras.

Using the formulas \eqref{der} and \eqref{formula} we have the following result.

\begin{pro}\label{derivadas}
The derivatives of the elements from the list \eqref{base2} are
given in the Figure 2.

\begin{figure}[h]
\begin{center}
\begin{tabular}{|c|c|c|}
  \hline
   & $\partial_1$ & $\partial_2$ \\
  \hline
  $[A]$ & $1$ & $0$ \\
  $[B]$ & $(1-q^{-3})x_2$ & $0$ \\
  $[C]$ & $q^2(1-q^{-3})^2x_2[D]+p_{21}(1-q^{-3})(q^3-q^2-q)[E]$ & $0$ \\
  $[D]$ & $(1-q^{-3})(1-q^{-2})x_2^2$ & $0$ \\
  $[E]$ & $(1-q^{-3})(1-q^{-2})(1-q^{-1})x_2^3$ & $0$ \\
  $[F]$ & $0$ & $1$ \\
  \hline
\end{tabular}
\caption{Table of Derivatives}
\end{center}
\end{figure}
\end{pro}

\emph{Proof.} Since $[A]=x_1$ and $[F]=x_2$, from the definition,
$\partial_1([A])=1$, $\partial_2([A])=0$, $\partial_1([F])=0$,
$\partial_2([F])=1$.

For $[B]$ we have $[B]=[x_1,x_2]=x_1x_2-p_{12}x_2x_1$. Using
(\ref{der}),
\begin{align*}
\partial_1([B])&=\partial_1(x_1x_2)-p_{12}\partial_1(x_2x_1)=\\&=\partial_1(x_1)x_2+p_{11}x_1\partial_1(x_2)-p_{12}(\partial_1(x_2)x_1+p_{21}x_2\partial_1(x_1))=\\&=x_2-p_{12}p_{21}x_2=(1-q^{-3})x_2.
\end{align*}
For $[D]=[[x_1,x_2],x_2]=[[B],x_2]=[B]x_2-p_{12}p_{22}x_2[B]$, we
have
\begin{align*}
\partial_1([D])&=\partial_1([B]x_2)-p_{12}p_{22}\partial_1(x_2[B])=\\&=\partial_1([B])x_2+p_{11}p_{21}[B]\partial_1(x_2)-p_{12}p_{22}(\partial_1(x_2)[B]+p_{21}x_2\partial_1([B]))=\\&=(1-p_{12}p_{21})x_2^2+0-0-p_{12}p_{22}p_{21}(1-p_{12}p_{21})x_2^2=\\&=(1-p_{12}p_{22}p_{21})(1-p_{12}p_{21})x_2^2=(1-q^{-3})(1-q^{-2})x_2^2.
\end{align*}
Again, for
$[E]=[[[x_1,x_2],x_2],x_2]=[[D],x_2]=[D]x_2-p_{12}p_{22}^2x_2[D]$,
we have
\begin{align*}
\partial_1([E])&=\partial_1([D]x_2)-p_{12}p_{22}^2\partial_1(x_2[D])=\\&=\partial_1([D])x_2+p_{11}p_{21}^2[D]\partial_1(x_2)-p_{12}p_{22}^2(\partial_1(x_2)[D]+p_{21}x_2\partial_1([D]))=\\&=(1-p_{12}p_{22}^2p_{21})(1-p_{12}p_{22}p_{21})(1-p_{12}p_{21})x_2^3=\\&=(1-q^{-3})(1-q^{-2})(1-q^{-1})x_2^3.
\end{align*}
Finally, for
$[C]=[[x_1,x_2],[[x_1,x_2],x_2]]=[B][D]-p_{11}p_{12}^2p_{21}p_{22}^2[D][B]$,
we have
\begin{align*}
\partial_1([C])&=\partial_1([B][D])-p_{11}p_{12}^2p_{21}p_{22}^2\partial_1([D][B])=\\&=\partial_1([B])[D]+p_{11}p_{21}[B]\partial_1([D])-\\&-p_{11}p_{12}^2p_{21}p_{22}^2(\partial_1([D])[B]+p_{11}p_{21}^2[D]\partial_1([B]))=\\&=(1-p_{12}p_{21})x_2[D]+p_{11}p_{21}(1-p_{12}p_{22}p_{21})(1-p_{12}p_{21})[B]x_2^2-\\&-p_{12}p_{22}^2(1-p_{12}p_{22}p_{21})(1-p_{12}p_{21})x_2^2[B]-p_{21}p_{22}^2(1-p_{12}p_{21})[D]x_2.
\end{align*}
Note that the elements $[D]x_2$ and $[B]x_2^2$ are not basis
elements, since $[B]>x_2$ and $[D]>x_2$. To write
$\partial_1([C])$ in the PBW-basis, we use that
\[[E]=[D]x_2-p_{12}p_{22}^2x_2[D],\]
\[[D]=[B]x_2-p_{12}p_{22}x_2[B],\] what provides
\[[D]x_2=[E]+p_{12}p_{22}^2x_2[D],\]
\[[B]x_2^2=[D]x_2+p_{12}p_{22}x_2[B]x_2=[E]+p_{12}p_{22}(1+p_{22})x_2[D]+p_{12}^2p_{22}^2x_2^2[B].\]
Using these relations we finally obtain
\begin{align*}
\partial_1([C])&=(1-p_{12}p_{21})x_2[D]+p_{11}p_{21}(1-p_{12}p_{22}p_{21})(1-p_{12}p_{21})([E]+\\&+p_{12}p_{22}(1+p_{22})x_2[D]
+p_{12}^2p_{22}^2x_2^2[B])-\\&-p_{12}p_{22}^2(1-p_{12}p_{22}p_{21})(1-p_{12}p_{21})x_2^2[B]
-p_{21}p_{22}^2(1-p_{12}p_{21})([E]+\\&+p_{12}p_{22}^2x_2[D])=(1-p_{12}p_{21})x_2[D](1
+p_{22}(1+p_{22})(1-p_{12}p_{22}p_{21})-\\&-p_{22})+p_{21}(1-p_{12}p_{21})[E](p_{11}(1-p_{12}p_{22}p_{21})-p_{22}^2)
+\\&+p_{12}p_{22}^2(1-p_{12}p_{22}p_{21})(1-p_{12}p_{21})x_2^2[B](p_{12}p_{11}p_{21}-1)
=\\&=q^2(1-q^{-3})^2x_2[D]+p_{21}(1-q^{-3})(q^3-q^2-q)[E].
\end{align*}

To calculate $\partial_2$, we note that directly from formula \eqref{formula} and $\partial_2(x_1)=0$ we have:
$$\partial_2([B])=\partial_2([D])=\partial_2([E])=0.$$
Now, for $[C]$:
\begin{align*}
\partial_2([C])&=\partial_2([B][D])-p_{11}p_{12}^2p_{21}p_{22}^2\partial_2([D][B])=\\&=\partial_2([B])[D]+p_{12}p_{22}[B]\partial_2([D])
-p_{12}p_{22}^2(\partial_2([D])[B]+\\&+p_{12}p_{22}^2[D]\partial_2([B]))=0.
\end{align*}
\hfill $\square$ \vspace{0.5cm}

\begin{obs}\label{puu}
It is proved in \cite[Proposition 4.7]{A} that for $[A], [B], [C], [D], [E], [F]$ from list \eqref{base2}, we
have $p(A,A)=q^3$, $p(B,B)=q$, $p(C,C)=q^3$, $p(D,D)=q$,
$p(E,E)=q^3$ and $p(F,F)=q$. However, this can easily be obtained with a simple calculation.
\end{obs}

\begin{teo}\label{pbw1}
If $q$ is not a root of $1$, then the values in
$U_q^+(\mathfrak{g})$ of the super-letters
\begin{align}\label{base2}
[A]&=x_1,\notag\\ [B]&=[x_1,x_2],\notag\\
[C]&=[[x_1,x_2],[[x_1,x_2],x_2]],\notag\\ [D]&=[[x_1,x_2],x_2],\notag\\
[E]&=[[[x_1,x_2],x_2],x_2],\notag\\ [F]&=x_2.\notag
\end{align}
form a set of PBW-generators for $U_q^+(\mathfrak{g})$ over
\emph{$\textbf{k}[G]$}, and each super-letter has infinite height.
If we suppose that $x_1>x_2$, then $A > B > C > D
> E > F$.
\end{teo}

\emph{Proof.} This
statement follows from \cite[Section 4.2]{A} due to the fact that $U_q^+({\mathfrak
g})$ is a  bosonization of a Nichols algebra generated by $x_1,x_2$. It is also a consequence of Proposition \ref{base},
since all the
hard in $U_q^+(\mathfrak{g})$ super-letters are contained in the list \eqref{base2}. If they are all hard and not zero, from Theorem
\ref{hard}, they form a set of PBW-generators for
$U_q^+(\mathfrak{g})$ over $\textbf{k}[G]$. Now we only have to
see that all heights are infinite.

From Remark \ref{puu}, we know that for every hard super-letter
$[u]$, either $p(u,u)=q$ or $p(u,u)=q^3$. But we are supposing
that $q$ is not a root of $1$, so $p(u,u)$ is not a primitive
$t$-th root of $1$ for any $t$, and from Definition \ref{altura}
we have that $h([u])$ is infinite. \hfill $\square$ \vspace{0.5cm}

Before we go to the next lemma, let us make some considerations
that will be necessary to prove it.

A super-letter $[u]$ is hard in $H$ if its value in $H$ is not a
linear combination of super-words of the same degree (\ref{grau})
in super-letters smaller than $[u]$. Let us notice that it
suffices to check that they are not a linear combination of
super-words of the same degree in hard super-letters smaller than $[u]$. Suppose that a non-hard super-letter $[u]$ is a
linear combination
$$[u]=\sum_i \alpha_i[v_{i_1}][v_{i_2}]\ldots[v_{i_k}]$$
where $[v_{i_j}]<[u]$ for every $i,j$ and the degree of
$[v_{i_1}][v_{i_2}]\ldots[v_{i_k}]$ is the same as the degree of
$[u]$. If one of the super-letters $[v_{i_j}]$ is not hard, we may
substitute
$$[v_{i_j}]=\sum_l \beta_l[w_{l_1}][w_{l_2}]\ldots[w_{l_m}]$$
where $[w_{l_n}]<[v_{i_j}]$ for every $l,n$ and the degree of
$[w_{l_1}][w_{l_2}]\ldots[w_{l_m}]$ is the same as the degree of
$[v_{i_j}]$. So we have $[w_{l_n}]<[v_{i_j}]<[u]$ and the degree
of $[w_{l_n}]$ is less than or equal to the degree of $[v_{i_j}]$
that is less than or equal to the degree of $[u]$. We have only a
finite number of super-letters smaller than $[u]$ with degree less
than or equal to the degree of $[u]$, so this process has to stop.
It means we can suppose that all $[v_{i_j}]$ are hard in
$U_q^+(\mathfrak{g})$, after making the necessary substitutions.

We know from Proposition \ref{base} that all hard in
$u_q^+(\mathfrak{g})$ super-letters belong to $\{[A], [B], [C],
[D], [E], [F]\}$. Since the super-letter $[[C],[F]]$ is not hard, it is a linear
combination of super-words of the same degree in hard super-letters smaller than
$[[C],[F]]$, which are $[D]$, $[E]$ and $[F]$.
The degree of $[[C],[F]]$ is $(2,4)$. The only possible
combination with $[D]$, $[E]$ and $[F]$ that has degree $(2,4)$ is
$[D]^2$. We conclude that
\begin{equation}\label{cf}
[[C],[F]]=\alpha [D]^2, \quad \alpha \in \textbf{k}.
\end{equation} In the same way, $[[C],[E]]$ is not hard and has
degree $(3,6)$, that provides
\begin{equation}\label{ce}
[[C],[E]]=\beta [D]^3, \quad \beta \in \textbf{k}. \end{equation}
Now let us see in the same way that
\begin{equation}\label{bc}
[[B],[C]]=[[C],[D]]=[[D],[E]]=0. \end{equation} The super-letter
$[[B],[C]]$ is not hard and has degree $(3,4)$. But there is no
combination of $[C]$, $[D]$, $[E]$ and $[F]$ (the hard
super-letters smaller than $[[B],[C]]$) that has this degree. So,
we have that $[[B],[C]]=0$. The same method can be used for
$[[C],[D]]=[[D],[E]]=0$.

\begin{lem}\label{dercol}
Let $[u]$ be an element from list \eqref{base2}. We have
\begin{equation*}
\underbrace{[[u],[[u], \ldots[[u]}_{l},\partial_i([u])]\ldots]]=0,
\end{equation*} for $l=1$ if $[u] \in \{[A], [E], [F]\}$, $l=2$ if $[u]=[C]$, and $l=3$ if $[u] \in \{[B], [D]\}$.
\end{lem}

\emph{Proof.} First we consider $[u]=[A]=x_1$. We have
$[[A],\partial_1([A])]=[x_1,1]=0$.

If $[u]=[B]$, then
$$[[B],\partial_1([B])]=(1-q^{-3})[[B],x_2]=(1-q^{-3})[D],$$ $$[[B],[D]]=[C].$$
Since from \eqref{bc} we have $[[B],[C]]=0$, we obtain the
required equality for $l=3$.

In the case $[u]=[C]$, using \eqref{cf}, \eqref{bc}, \eqref{ce}
and \eqref{bracprod2} we have
$$[[C],\partial_1([C])]=\alpha[[C],x_2]\cdot[D]+\beta
x_2\cdot[[C],[D]]+\gamma[[C],[E]]=\delta[D]^3,$$
$$[[C],[D]^3]=\varepsilon[[C],[D]]\cdot[D]^2+\theta[D]\cdot[[C],[D]]\cdot[D]+\lambda[D]^2\cdot[[C],[D]]=0,$$
 with $\alpha, \beta, \gamma, \delta, \varepsilon, \theta, \lambda \in
 \textbf{k}$.

If $[u]=[D]$, then
$$[[D],x_2]=[E],\quad [[D],[E]]=0,$$
so we obtain from \eqref{bracprod2} the following relations
\begin{align*}
[[D],[[D],x_2^2]]&=[[D],[[D],x_2]\cdot x_2]+[[D],\alpha
x_2\cdot[[D],x_2]]=\\&=[[D],[E]x_2]+[[D],\alpha
x_2[E]]=\\&=[E]\cdot[[D],x_2]+\beta[[D],[E]]\cdot x_2+\alpha
x_2\cdot[[D],[E]]+\gamma[[D],x_2]\cdot[E]=\\&=\delta[E]^2,
\end{align*}
with $\alpha, \beta, \gamma, \delta \in \textbf{k}$ and
$$[[D],[E]^2]=[E]\cdot[[D],[E]]+\varepsilon[[D],[E]]\cdot[E]=0,$$
$\varepsilon \in \textbf{k}$, that gives us
$$[[D],[[D],[[D],\partial_1([D])]]]=0.$$

For $[u]=[E]$ we have
$$[[E],x_2]=0,$$
so from \eqref{bracprod2} it follows
$$[[E],x_2^3]=[[E],x_2]\cdot x_2^2+\alpha x_2\cdot[[E],x_2]\cdot x_2+\beta
x_2^2\cdot[[E],x_2]=0,$$ $\alpha,\beta \in \textbf{k}$. Hence
$[[E],\partial_1([E])]=0$.

Our last possibility is $[u]=[F]=x_2$. In this case,
$$[[F],\partial_2([F])]=[x_2,1]=0.$$

It only remains to see that
$$[[u],\partial_2([u])]=0$$
for $[u] \in \{[A],[B],[C],[D],[E]\}$ and
$$[[F],\partial_1([F])]=0.$$ This is obvious, since we have
$\partial_2([A])=\partial_2([B])=\partial_2([C])=\partial_2([D])=\partial_2([E])=\partial_1([F])=0$.
\hfill $\square$ \vspace{0.5cm}

\begin{teo}\label{pbwq}
If $q$ has finite multiplicative order $t$, $t>4$, $t\neq 6$, then
the values in $u_q^+(\mathfrak{g})$ of the super-letters from list
\eqref{base2} form a set of PBW-generators for
$u_q^+(\mathfrak{g})$ over \emph{$\textbf{k}[G]$}. The height $h$
of $[u] \in \{[B], [D], [F]\}$ equals $t$. For $[u] \in \{[A],
[C], [E]\}$ we have $h=t$ if $3$ is not a divisor of $t$ and
$h=\frac{t}{3}$ otherwise. In all cases $[u]^h=0$ in
$u_q^+(\mathfrak{g})$.
\end{teo}

\emph{Proof.} In the same way as in Theorem \ref{pbw1} we see that the
super-letters from the list \eqref{base2} form a set of PBW-generators for
$u_q^+(\mathfrak{g})$ over $\textbf{k}[G]$. Now
let us examine their heights.

First we notice that, if $p(u,u)$ is a primitive $t_u$-th root of
$1$ and
\begin{equation*}
\underbrace{[u,[u, \ldots[u}_{t_u-1},\partial_i(u)]\ldots]]=0,
\end{equation*} then from Lemma \ref{derivadasq1} we have
$\partial_i([u]^{t_u})=0$ in $u_q^+(\mathfrak{g})$.

In the case $[u] \in \{[B],[D],[F]\}$, we have $p(u,u)=q$. So
$t_u=t$, since $q$ is a primitive $t$-th root of $1$. From Lemmas
\ref{derivadasq1} and \ref{dercol}, we have
$\partial_i([u]^{t})=0$ in $u_q^+(\mathfrak{g})$ for $i=1,2$ and
$t\geq 5$. Now we apply the Milinski-Schneider criterion (Lemma
\ref{MS}), and obtain $[u]^{t}=0$. We get that $t$ is the height
of $[B],[D],[F]$.

For $[u] \in \{[A],[C],[E]\}$, we have $p(u,u)=q^3$. Again $q$ is
a primitive $t$-th root of $1$, so $t_u$ equals $t$ if $3$ is not
a divisor of $t$ and $\frac{t}{3}$ otherwise. From Lemmas
\ref{derivadasq1} and \ref{dercol}, we have
$\partial_i([u]^{t_u})=0$ in $u_q^+(\mathfrak{g})$ for $i=1,2$ and
$t_u\geq 3$. Since $t\geq 5$ and $t\neq 6$, its suffices to have
$t_u\geq 3$. Again, the Milinski-Schneider criterion provides
$[u]^{t_u}=0$. So the height of $[A],[C],[E]$ is $t$ or
$\frac{t}{3}$. \hfill $\square$ \vspace{0.5cm}

\begin{cor}\label{s1}
The exponent $s$ given in \eqref{geradores} is $1$ for every
PBW-generator $[u]$.
\end{cor}

\emph{Proof.} From Lemma \ref{s1t}, we know that either $s=1$ or $p(u,u)$
is a primitive $t$-th root of $1$ and $s=t$, or (in the case of
positive characteristic) $s=t(char \textbf{k})^r$. Since from
Theorem \ref{pbwq} we have $[u]^t=0$ if $p(u,u)$ is a primitive
$t$-th root of $1$, we obtain $s=1$. \hfill $\square$
\vspace{0.5cm}

\section{Lattice of Right Coideal Subalgebras}

\quad In this section we are going to describe all the
(homogeneous) right coideal subalgebras containing $\textbf{k}[G]$
of the multiparameter quantum group $U_q^+(\mathfrak{g})$
(respectively, of $u_q^+(\mathfrak{g})$), where $G$ is the set of
group-like elements and $\mathfrak{g}$ is the simple Lie algebra
of type $G_2$.

\vspace{0.5cm}
\emph{Proof of Theorem \ref{lattice}.} Let $\textbf{U}$ be a right coideal subalgebra of
$U_q^+(\mathfrak{g})$ (respectively, $u_q^+(\mathfrak{g})$). From
Proposition \ref{pbwU} and Corollary \ref{s1}, the PBW-generators
for $\textbf{U}$ have the form
$$[u]+\sum \alpha_i W_i$$
where $[u]$ is a hard super-letter, $W_i$ are the basis
super-words starting with super-letters smaller than $[u]$, and
$D(W_i)=D([u])$. So, all the possibilities are:

$x_1,$

$[B]+\alpha x_2x_1,$

$[D]+\alpha x_2^2x_1+\beta x_2[B],$

$[E]+\alpha x_2^3x_1+\beta x_2^2[B]+\gamma x_2[D],$

$[C]+\alpha x_2^3x_1^2+\beta x_2^2[B]x_1+\gamma x_2[D]x_1+\delta
x_2[B]^2+\varepsilon [D][B]+\tau [E]x_1,$

$x_2.$

Moreover, according to the construction, the set of PBW-generators
has not more than one PBW-generator of each of the six mentioned
types. In particular, each proper right coideal subalgebra of
$U_q^+(\mathfrak{g})$ or $u_q^+(\mathfrak{g})$ has not more than
five PBW-generators.

Our first goal is to calculate all the possible values for the
coefficients above.

Suppose that $[B]+\alpha x_2x_1$ is a generator of $\textbf{U}$.
Since $U_A=\textbf{U}\cap A$ is a differential subalgebra, the
following elements are in $\textbf{U}$:
\begin{enumerate}
    \item $\partial_2([B]+\alpha x_2x_1)=\alpha x_1,$
    \item $\partial_1([B]+\alpha x_2x_1)=(1-q^{-3}+\alpha p_{21})x_2.$
\end{enumerate}

If both $\alpha \neq 0$ and $(1-q^{-3}+\alpha p_{21}) \neq 0$,
then $x_1$ and $x_2$ are in $\textbf{U}$, and
$\textbf{U}=U_q^+(\mathfrak{g})$. So we may suppose that $\alpha
=0$ or $(1-q^{-3}+\alpha p_{21})=0$. In the first case, $[B]$ is a
generator for $\textbf{U}$. In the second case, $\alpha
=(q^{-3}-1)/p_{21}$ and the generator is
$$[B]+\alpha x_2x_1=[B]+\frac{(q^{-3}-1)x_2x_1}{p_{21}}=x_1x_2-p_{21}^{-1}x_2x_1=-p_{21}^{-1}[x_2,x_1].$$

Suppose now that $[D]+\alpha x_2^2x_1+\beta x_2[B]$ is a generator
of $\textbf{U}$. Then $\textbf{U}$ has the following elements:
\begin{enumerate}
    \item $\partial_2([D]+\alpha x_2^2x_1+\beta x_2[B])=\alpha
(1+q)x_2x_1+\beta[B],$
    \item $\partial_2^2([D]+\alpha x_2^2x_1+\beta x_2[B])=\alpha (1+q)x_1,$
    \item $\partial_1 ([D]+\alpha x_2^2x_1+\beta
x_2[B])=((1-q^{-3})(1-q^{-2})+\alpha p_{21}^2+\beta
p_{21}(1-q^{-3}))x_2^2,$
    \item $\partial_1\partial_2([D]+\alpha x_2^2x_1+\beta x_2[B])=(\alpha
(1+q)p_{21}+\beta(1-q^{-3}))x_2.$
\end{enumerate}

Again we have two possibilities. One is that $x_2 \in \textbf{U}$.
In this case, from the second line, $\alpha=0$. If $\beta=0$, then
$[D]$ is a generator. If $\beta\neq 0$, then $[D]$ may also be
considered as a generator. Thus $\partial_2([D]+\beta
x_2[B])=\beta[B]$, so $[B],x_2 \in \textbf{U}$ imply $[D] \in
\textbf{U}$. The second possibility is that $x_1 \in \textbf{U}$,
and from lines $3$ and $4$ we get
\begin{itemize}
\item $(1-q^{-3})(1-q^{-2})+\alpha p_{21}^2+\beta
p_{21}(1-q^{-3})=0,$ \item $\alpha (1+q)p_{21}+\beta(1-q^{-3})=0.$
\end{itemize}

Solving this system (of two equations and two variables) we find
$$\alpha=\frac{(1-q^{-3})(1-q^{-2})}{qp_{21}^2},\quad \beta=-\frac{(1-q^{-2})(1+q)}{p_{21}q}$$
and the generator is
$$[D]+\frac{(1-q^{-3})(1-q^{-2})}{qp_{21}^2} x_2^2x_1-\frac{(1-q^{-2})(1+q)}{p_{21}q} x_2[B]=q^{-1}p_{21}^{-2}[x_2,[x_2,x_1]].$$

If $[E]+\alpha x_2^3x_1+\beta x_2^2[B]+\gamma x_2[D]$ is a
generator of $\textbf{U}$, we have:
\begin{enumerate}
\item $\partial_2([E]+\alpha x_2^3x_1+\beta x_2^2[B]+\gamma
x_2[D])=\alpha(1+q+q^2)x_2^2x_1+\beta(1+q)x_2[B]+\gamma[D],$

\item $\partial_2^2([E]+\alpha x_2^3x_1+\beta x_2^2[B]+\gamma
x_2[D])=\alpha(1+q+q^2)(1+q)x_2x_1+\beta(1+q)[B],$

\item $\partial_2^3([E]+\alpha x_2^3x_1+\beta x_2^2[B]+\gamma
x_2[D])=\alpha(1+q+q^2)(1+q)x_1,$

\item $\partial_1\partial_2^2([E]+\alpha x_2^3x_1+\beta
x_2^2[B]+\gamma
x_2[D])=(1+q)(\alpha(1+q+q^2)p_{21}+\beta(1-q^{-3}))x_2,$

\item $\partial_1\partial_2([E]+\alpha x_2^3x_1+\beta
x_2^2[B]+\gamma x_2[D])=(\alpha p_{21}^2(1+q+q^2)+\beta
p_{21}(1+q)(1-q^{-3})+\gamma(1-q^{-3})(1-q^{-2}))x_2^2,$

\item $\partial_1([E]+\alpha x_2^3x_1+\beta x_2^2[B]+\gamma
x_2[D])=((1-q^{-3})(1-q^{-2})(1-q^{-1})+\alpha p_{21}^3+\beta
p_{21}^2(1-q^{-3})+\gamma p_{21}(1-q^{-3})(1-q^{-2}))x_2^3.$
\end{enumerate}

If $x_2 \in \textbf{U}$, from line $3$ we have $\alpha=0$ and
three possibilities remain. If $\beta=\gamma=0$, then $[E]$ is a
generator. If $\beta=0$ and $\gamma \neq 0$, then $[E]+\gamma
x_2[D] \in \textbf{U}$, while
$$\partial_1([E]+\gamma x_2[D])=(1-q^{-3})(1-q^{-2})((1-q^{-1})+\gamma p_{21})x_2^3,$$
$$\partial_2([E]+\gamma x_2[D])=\gamma[D].$$
Hence $[E] \in \textbf{U}$ and we still may consider the generator
to be $[E]$. If $\beta \neq 0$, then $[E]+\beta x_2^2[B]+\gamma
x_2[D] \in \textbf{U}$. In this case
$$\partial_2^2([E]+\beta
x_2^2[B]+\gamma x_2[D])=\beta(1+q)[B]$$
$$\partial_1([B])=(1-q^{-3})x_2$$
provide $[E] \in \textbf{U}$.

If $x_1 \in \textbf{U}$, then from lines $4$, $5$ and $6$ we get
\begin{itemize}
\item $\alpha(1+q+q^2)p_{21}+\beta(1-q^{-3})=0,$ \item $\alpha
p_{21}^2(1+q+q^2)+\beta
p_{21}(1+q)(1-q^{-3})+\gamma(1-q^{-3})(1-q^{-2})=0,$ \item
$(1-q^{-3})(1-q^{-2})(1-q^{-1})+\alpha p_{21}^3+\beta
p_{21}^2(1-q^{-3})+\gamma p_{21}(1-q^{-3})(1-q^{-2})=0.$
\end{itemize}

Solving this system we have that the generator $[E]+\alpha
x_2^3x_1+\beta x_2^2[B]+\gamma x_2[D]$ is a multiple of
$[x_2,[x_2,[x_2,x_1]]]$.

The last possibility is that $[C]+\alpha x_2^3x_1^2+\beta
x_2^2[B]x_1+\gamma x_2[D]x_1+\delta x_2[B]^2+\varepsilon
[D][B]+\tau [E]x_1$ is a generator of $\textbf{U}$. Now we have

\begin{enumerate}

\item $\partial_2\partial_1\partial_2^2([C]+\alpha
x_2^3x_1^2+\beta x_2^2[B]x_1+\gamma x_2[D]x_1+\delta
x_2[B]^2+\varepsilon [D][B]+\tau [E]x_1)=(1+q)(\alpha
p_{21}(1+q+q^2)(1+q^3)+\beta(1+q^{-3}))x_1,$

\item $\partial_2^3\partial_1([C]+\alpha x_2^3x_1^2+\beta
x_2^2[B]x_1+\gamma x_2[D]x_1+\delta x_2[B]^2+\varepsilon
[D][B]+\tau [E]x_1)=(\alpha p_{21}^3(1+q^3)+\beta
p_{21}^2(1-q^{-3})+\gamma
p_{21}(1-q^{-3})(1-q^{-2})+\tau(1-q^{-3})(1-q^{-2})(1-q^{-1}))(1+q+q^2)(1+q)x_1,$

\item $\partial_2^3([C]+\alpha x_2^3x_1^2+\beta x_2^2[B]x_1+\gamma
x_2[D]x_1+\delta x_2[B]^2+\varepsilon [D][B]+\tau
[E]x_1)=\alpha(1+q+q^2)(1+q)x_1^2,$

\item $\partial_2^2\partial_1\partial_2([C]+\alpha
x_2^3x_1^2+\beta x_2^2[B]x_1+\gamma x_2[D]x_1+\delta
x_2[B]^2+\varepsilon [D][B]+\tau [E]x_1)=(\alpha
p_{21}^2(1+q+q^2)(1+q^3)+\beta
p_{21}(1+q)(1-q^{-3})+\gamma(1-q^{-3})(1-q^{-2}))(1+q)x_1,$

\item $\partial_1^2([C]+\alpha x_2^3x_1^2+\beta x_2^2[B]x_1+\gamma
x_2[D]x_1+\delta x_2[B]^2+\varepsilon [D][B]+\tau
[E]x_1)=p_{21}(p_{21}^2(\alpha p_{21}^3(1+q^3)+\beta
p_{21}^2(1-q^{-3})+\gamma
p_{21}(1-q^{-3})(1-q^{-2})+\tau(1-q^{-3})(1-q^{-2})(1-q^{-1}))+p_{21}(1-q^{-3})(\beta
p_{21}^3q^3+\delta
p_{21}(1-q^{-3})+\varepsilon(1-q^{-3})(1-q^{-2})+\delta
p_{21}q(1-q^{-3}))+(1-q^{-3})(1-q^{-2})(\gamma
p_{21}^3q^3+q^2(1-q^{-3})^2+\delta
p_{21}^2q^3(1-q^{-3})+\varepsilon
p_{21}q^2(1-q^{-3}))+(1-q^{-3})(1-q^{-2})(1-q^{-1})q^3(\varepsilon
p_{21}(1-q^{-3})+\tau
p_{21}^2+(1-q^{-3})(1-q^{-1}-q^{-2})))x_2^3,$

\item $\partial_1\partial_2\partial_1([C]+\alpha x_2^3x_1^2+\beta
x_2^2[B]x_1+\gamma x_2[D]x_1+\delta x_2[B]^2+\varepsilon
[D][B]+\tau [E]x_1)=(p_{21}^2(1+q+q^2)(\alpha
p_{21}^3(1+q^3)+\beta p_{21}^2(1-q^{-3})+\gamma
p_{21}(1-q^{-3})(1-q^{-2})+\tau(1-q^{-3})(1-q^{-2})(1-q^{-1}))+p_{21}(1+q)(1-q^{-3})(\beta
p_{21}^3q^3+\delta
p_{21}(1-q^{-3})+\varepsilon(1-q^{-3})(1-q^{-2})+\delta
p_{21}q(1-q^{-3}))+(1-q^{-3})(1-q^{-2})q^2(\gamma
p_{21}^3q+(1-q^{-3})^2+\delta p_{21}^2q(1-q^{-3})+\varepsilon
p_{21}(1-q^{-3})))x_2^2,$

\item $\partial_1\partial_2^2\partial_1([C]+\alpha
x_2^3x_1^2+\beta x_2^2[B]x_1+\gamma x_2[D]x_1+\delta
x_2[B]^2+\varepsilon [D][B]+\tau [E]x_1)=(p_{21}(1+q+q^2)(\alpha
p_{21}^3(1+q^3)+\beta p_{21}^2(1-q^{-3})+\gamma
p_{21}(1-q^{-3})(1-q^{-2})+\tau(1-q^{-3})(1-q^{-2})(1-q^{-1}))+(1-q^{-3})(\beta
p_{21}^3q^3+\delta
p_{21}(1-q^{-3})+\varepsilon(1-q^{-3})(1-q^{-2})+\delta
p_{21}q(1-q^{-3})))(1+q)x_2,$

\item $\partial_1^2\partial_2^2([C]+\alpha x_2^3x_1^2+\beta
x_2^2[B]x_1+\gamma x_2[D]x_1+\delta x_2[B]^2+\varepsilon
[D][B]+\tau [E]x_1)=p_{21}(1+q)(1+q^3)(\alpha
p_{21}(1+q+q^2)+\beta(1-q^{-3}))x_2,$

\item $\partial_1^2\partial_2([C]+\alpha x_2^3x_1^2+\beta
x_2^2[B]x_1+\gamma x_2[D]x_1+\delta x_2[B]^2+\varepsilon
[D][B]+\tau [E]x_1)=p_{21}(p_{21}(\alpha
p_{21}^2(1+q+q^2)(1+q^3)+\beta
p_{21}(1+q)(1-q^{-3})+\gamma(1-q^{-3})(1-q^{-2}))+(1-q^{-3})(1+q)(\beta
p_{21}^2q^3+\delta(1-q^{-3}))+(1-q^{-3})(1-q^{-2})q^3(\gamma
p_{21}+\delta(1-q^{-3})))x_2^2,$

\item $\partial_1\partial_2\partial_1\partial_2([C]+\alpha
x_2^3x_1^2+\beta x_2^2[B]x_1+\gamma x_2[D]x_1+\delta
x_2[B]^2+\varepsilon [D][B]+\tau [E]x_1)=(1+q)(p_{21}(\alpha
p_{21}^2(1+q+q^2)(1+q^3)+\beta
p_{21}(1+q)(1-q^{-3})+\gamma(1-q^{-3})(1-q^{-2}))+(1-q^{-3})(\beta
p_{21}^2q^3+\delta(1-q^{-3})))x_2.$

\end{enumerate}

If $x_2 \in \textbf{U}$, then from the first $4$ equalities we
have
\begin{itemize}
\item $\alpha p_{21}(1+q+q^2)(1+q^3)+\beta(1+q^{-3})=0,$ \item
$\alpha p_{21}^3(1+q^3)+\beta p_{21}^2(1-q^{-3})+\gamma
p_{21}(1-q^{-3})(1-q^{-2})+\tau(1-q^{-3})(1-q^{-2})(1-q^{-1})=0,$
\item $\alpha=0,$ \item $\alpha p_{21}^2(1+q+q^2)(1+q^3)+\beta
p_{21}(1+q)(1-q^{-3})+\gamma(1-q^{-3})(1-q^{-2})=0.$
\end{itemize}

Solving this system we find $\alpha=\beta=\gamma=\tau=0$. In this
case the generator is $[C]$. If $\delta=\varepsilon=0$, it is
obvious. If one of $\delta,\varepsilon$ is not zero, then
$$\partial_2^2\partial_1([C]+\delta x_2[B]^2+\varepsilon[D][B])=(1-q^{-3})(1+q)(\delta
    p_{21}(1+q)+\varepsilon(1-q^{-2}))[B]$$
$$\partial_1([B])=(1-q^{-3})x_2$$
$$[D]=[[B],x_2]$$
imply that $[C]$ still belongs to $\textbf{U}$. Notice that if
$\delta p_{21}(1+q)+\varepsilon(1-q^{-2})=0$, then we have
$\delta=-\frac{\varepsilon(1-q^{-2})}{p_{21}(1+q)}$. Therefore,
$\delta \neq 0$ and $\varepsilon \neq 0$, and since
$\partial_2\partial_1\partial_2([C]+\delta
x_2[B]^2+\varepsilon[D][B])=-\varepsilon(1-q^{-3})(1-q^{-2})p_{21}^{-1}[B]$,
we obtain $[B] \in \textbf{U}$.

If $x_1 \in \textbf{U}$, from equalities $5-10$ we have
\begin{itemize}
\item $p_{21}^2(\alpha p_{21}^3(1+q^3)+\beta
p_{21}^2(1-q^{-3})+\gamma
p_{21}(1-q^{-3})(1-q^{-2})+\tau(1-q^{-3})(1-q^{-2})(1-q^{-1}))+p_{21}(1-q^{-3})(\beta
p_{21}^3q^3+\delta
p_{21}(1-q^{-3})+\varepsilon(1-q^{-3})(1-q^{-2})+\delta
p_{21}q(1-q^{-3}))+(1-q^{-3})(1-q^{-2})q^2(\gamma
p_{21}^3q+(1-q^{-3})+\delta p_{21}^2q(1-q^{-3})+\varepsilon
p_{21}(1-q^{-3}))+(1-q^{-3})(1-q^{-2})(1-q^{-1})q^3(\varepsilon
p_{21}(1-q^{-3})+\tau p_{21}^2+(1-q^{-3})(1-q^{-1}-q^{-2}))=0,$

\item $p_{21}^2(1+q+q^2)(\alpha p_{21}^3(1+q^3)+\beta
p_{21}^2(1-q^{-3})+\gamma
p_{21}(1-q^{-3})(1-q^{-2})+\tau(1-q^{-3})(1-q^{-2})(1-q^{-1}))+p_{21}(1+q)(1-q^{-3})(\beta
p_{21}^3q^3+\delta
p_{21}(1-q^{-3})+\varepsilon(1-q^{-3})(1-q^{-2})+\delta
p_{21}q(1-q^{-3}))+(1-q^{-3})(1-q^{-2})q^2(\gamma
p_{21}^3q+(1-q^{-3})+\delta p_{21}^2q(1-q^{-3})+\varepsilon
p_{21}(1-q^{-3}))=0,$

\item $p_{21}(1+q+q^2)(\alpha p_{21}^3(1+q^3)+\beta
p_{21}^2(1-q^{-3})+\gamma
p_{21}(1-q^{-3})(1-q^{-2})+\tau(1-q^{-3})(1-q^{-2})(1-q^{-1}))+(1-q^{-3})(\beta
p_{21}^3q^3+\delta
p_{21}(1-q^{-3})+\varepsilon(1-q^{-3})(1-q^{-2})+\delta
p_{21}q(1-q^{-3}))=0,$

\item $\alpha p_{21}(1+q+q^2)+\beta (1-q^{-3})=0,$

\item $p_{21}(\alpha p_{21}^2(1+q+q^2)(1+q^3)+\beta
p_{21}(1+q)(1-q^{-3})+\gamma(1-q^{-3})(1-q^{-2}))+(1-q^{-3})(1+q)(\beta
p_{21}^2q^3+\delta(1-q^{-3}))+q^3(1-q^{-3})(1-q^{-2})(\gamma
p_{21}+\delta(1-q^{-3}))=0,$

\item $p_{21}(\alpha p_{21}^2(1+q+q^2)(1+q^3)+\beta
p_{21}(1+q)(1-q^{-3})+\gamma(1-q^{-3})(1-q^{-2}))+(1-q^{-3})(\beta
p_{21}^2q^3+\delta(1-q^{-3}))=0.$
\end{itemize}

Solving this system we find that the generator $[C]+\alpha
x_2^3x_1^2+\beta x_2^2[B]x_1+\gamma x_2[D]x_1+\delta
x_2[B]^2+\varepsilon [D][B]+\tau [E]x_1$ is a multiple of
$[[x_{1},x_{2}],[x_{2},[x_{2},x_{1}]]]$.

Now we are ready to study the possible right coideal subalgebras,
as we know that the possible PBW-generators for $\textbf{U}$ are

$x_1,$

$[B]$ or $[x_2,x_1]$,

$[D]$ or $[x_2,[x_2,x_1]]$,

$[E]$ or $[x_2,[x_2,[x_2,x_1]]]$,

$[C]$ or $[[x_1,x_2],[x_2,[x_2,x_1]]]$,

$x_2$.

Let us notice first that, by $\langle [u] \rangle$ we mean the
smallest right coideal subalgebra containing $[u]$ and
$\textbf{k}[G]$.

The right coideal subalgebra generated by $x_2$ has PBW-generators
$\{x_2\}$. In the same way, the right coideal subalgebra generated
by $x_1$ has PBW-generators $\{x_1\}$.

If $[E] \in \textbf{U}$, then $x_2 \in \textbf{U}$, since
$$\partial_1([E])=(1-q^{-3})(1-q^{-2})(1-q^{-2})x_2^3.$$
Hence $\langle [E] \rangle$ has PBW-generators $\{x_2, [E]\}.$

If $[D] \in \textbf{U}$, then $x_2$ and $[E]$ belong to
$\textbf{U}$, thus
$$\partial_1([D])=(1-q^{-3})(1-q^{-2})x_2^2,$$
$$[E]=[[D],x_2]$$
and $\langle [D] \rangle$ has PBW-generators $\{x_2, [D], [E]\}.$

If $[C] \in \textbf{U}$, then $x_2, [D]$ and $[E]$ are in
$\textbf{U}$, because
$$\partial_2 \partial_1([C])=q^2(1-q^{-3})^2[D],$$
providing that $\langle [C] \rangle$ has PBW-generators $\{x_2,
[C], [D], [E]\}.$

If $[B] \in \textbf{U}$, then $x_2, [C], [D], [E] \in \textbf{U}$,
since $$\partial_1([B])=(1-q^{-3})x_2,$$
$$[C]=[[B],[[B],x_2]]=[[B],[D]].$$
So, we obtain that the right coideal subalgebra generated by $[B]$
has PBW-generators $\{x_2, [B], [C], [D], [E]\}.$

If we include $x_1$ in any of these four right coideal
subalgebras, we have $\textbf{U}=U_q^+(\mathfrak{g})$, since $x_2$
belongs to all of them.

If $[x_2,x_1] \in \textbf{U}$, then $x_1 \in \textbf{U}$ since
\begin{align*}
\partial_2([x_2,x_1])&=\partial_2(x_2x_1)-p_{21}\partial_2(x_1x_2)=\\&=\partial_2(x_2)x_1+p_{22}x_2\partial_2(x_1)-p_{21}(\partial_2(x_1)x_2+p_{12}x_1\partial_2(x_2))=\\&=x_1-p_{21}p_{12}x_1=(1-q^{-3})x_1,
\end{align*}
and $\langle [x_2,x_1] \rangle$ has PBW-generators $\{x_1,
[x_2,x_1]\}$.

If $[x_2,[x_2,x_1]] \in \textbf{U}$, then $x_1, [x_2,x_1] \in
\textbf{U}$, thus
\begin{align*}
\partial_2([x_2,[x_2,x_1]])&=\partial_2(x_2[x_2,x_1])-p_{22}p_{21}\partial_2([x_2,x_1]x_2)=\\&=\partial_2(x_2)[x_2,x_1]+p_{22}x_2\partial_2([x_2,x_1])-\\&-p_{22}p_{21}(\partial_2([x_2,x_1])x_2+p_{22}p_{12}[x_2,x_1]\partial_2(x_2))=\\&=[x_2,x_1]+p_{22}(1-q^{-3})[x_2,x_1]-p_{21}p_{12}q^2[x_2,x_1]=\\&=(1+q)(1-q^{-2})[x_2,x_1].
\end{align*}
Hence $\langle [x_2,[x_2,x_1]] \rangle$ has PBW-generators $\{x_1,
[x_2,x_1], [x_2,[x_2,x_1]]\}.$

If $[x_2,[x_2,[x_2,x_1]]] \in \textbf{U}$, then $x_1, [x_2,x_1],
[x_2,[x_2,x_1]] \in \textbf{U}$ because
\begin{align*}
\partial_2([x_2,[x_2,[x_2,x_1]]])&=\partial_2(x_2[x_2,[x_2,x_1]])-p_{22}^2p_{21}\partial_2([x_2,[x_2,x_1]]x_2)=\\&=\partial_2(x_2)[x_2,[x_2,x_1]]+p_{22}x_2\partial_2([x_2,[x_2,x_1]])-\\&-p_{22}^2p_{21}(\partial_2([x_2,[x_2,x_1]])x_2+p_{22}^2p_{12}[x_2,[x_2,x_1]]\partial_2(x_2))=\\&=[x_2,[x_2,x_1]]+p_{22}(1+q)(1-q^{-2})[x_2,[x_2,x_1]]-\\&-p_{21}p_{12}q^4[x_2,[x_2,x_1]]=q^2(1-q^{-3})[x_2,[x_2,x_1]]
\end{align*}
implies that the PBW-generators are $\{x_1, [x_2,x_1],
[x_2,[x_2,x_1]], [x_2,[x_2,[x_2,x_1]]]\}.$

If $[[x_{1},x_{2}],[x_{2},[x_{2},x_{1}]]] \in \textbf{U}$, then
$x_1, [x_2,x_1], [x_2,[x_2,x_1]], [x_2,[x_2,[x_2,x_1]]] \in
\textbf{U}$ since
\begin{align*}
\partial_1([[x_{1},x_{2}],[x_{2},[x_{2},x_{1}]]])&=\partial_1([x_{1},x_{2}][x_{2},[x_{2},x_{1}]])-\\&-p_{12}^2p_{11}p_{22}^2p_{21}\partial_1([x_{2},[x_{2},x_{1}]][x_{1},x_{2}])=\\&=\partial_1([x_{1},x_{2}])[x_{2},[x_{2},x_{1}]]+\\&+p_{11}p_{21}[x_{1},x_{2}]\partial_1([x_{2},[x_{2},x_{1}]])-\\&-p_{12}q^2(\partial_1([x_{2},[x_{2},x_{1}]])[x_{1},x_{2}]+\\&+p_{21}^2p_{11}[x_{2},[x_{2},x_{1}]]\partial_1([x_{1},x_{2}]))=\\&=(1-q^{-3})x_2[x_{2},[x_{2},x_{1}]]+\\&+p_{21}^2p_{11}(1-q^{-3})[x_{2},[x_{2},x_{1}]]x_{2}=\\&=(1-q^{-3})[x_2,[x_{2},[x_{2},x_{1}]]].
\end{align*}
So $\{x_1, [x_2,x_1], [x_2,[x_2,x_1]], [x_2,[x_2,[x_2,x_1]]],
[[x_1,x_2],[x_2,[x_2,x_1]]]\}$ is a set of PBW-generators for
$\textbf{U}$.

Again, if we include $x_2$ in any of these four right coideal
subalgebras, we have $\textbf{U}=U_q^+(\mathfrak{g})$.

From the fact that a right coideal subalgebras can not have two
generators of the same type we conclude that these are all the
(homogeneous) right coideal subalgebras of $U_q^+(\mathfrak{g})$
(respectively, of $u_q^+(\mathfrak{g})$) that contain
$\textbf{k}[G]$, and we have the Figure 1. The theorem is proved.

\hfill $\square$ \vspace{0.5cm}

\section*{Acknowledgments}

I would like to thank professor Kharchenko for all his attention
during my research period in Mexico. I also thank him for
proposing this theme as part of my thesis and for all his suggestions to
this article.

\end{document}